\documentclass[12pt]{article}
\usepackage{latexsym}
\usepackage[centertags]{amsmath}
\usepackage{amsfonts}
\usepackage{amssymb}

      \newtheorem{theorem}{Theorem}[section]

      \newtheorem{lemma}[theorem]{Lemma}
      
      \newtheorem{definition}{Definition}[section]

      \def\nn{\nonumber}
      \def\rf#1{\mbox{$(\ref{#1})$}}

      \def\be{\begin{equation}} 
      \def\ee{\end{equation}} 
      \def\beqn{\begin{eqnarray}} 
      \def\eeqn{\end{eqnarray}} 
      \def\beq{\begin{eqnarray*}} 
      \def\eeq{\end{eqnarray*}}
      \def\proof{{\bf Proof:}\ }
      \def\mb{\mbox} 
      \def\ga{\gamma} 
      \def\ep{\varepsilon} 
      \def\la{\lambda} 
      \def\ra{\rightarrow} 

\newcommand{\re}{\mathrm{e}} 

\def\M{{\cal M}}

\def\ep{{\epsilon}}\def\ga{{\gamma}}

\def\<{\left<}\def\>{\right>}
\def\({\left(}\def\){\right)}

\newfam\msbmfam\font\tenmsbm=msbm10\textfont
\msbmfam=\tenmsbm\font\sevenmsbm=msbm7

\scriptfont\msbmfam=\sevenmsbm

\begin{document}

      \title{\bf Gamma-Dirichlet Structure and Two Classes of Measure-Valued Processes \thanks{Research supported by
      the Natural Science and Engineering Research Council of Canada}}
      \author{ Shui Feng and Fang Xu \\McMaster
      University}
      \date{\today}
      \maketitle
      \begin{abstract}
      The Gamma-Dirichlet structure corresponds to the decomposition of the gamma process into the independent product of a gamma random variable and a Dirichlet process. This structure allows us to study the properties of the Dirichlet process through the gamma process and vice versa. In this article, we begin with a brief review of existing results concerning the Gamma-Dirichlet structure.  New results are obtained for the large deviations of the jump sizes of the gamma process and the quasi-invariance of the two-parameter Poisson-Dirichlet distribution. The laws of the gamma process and the Dirichlet process are the respective reversible measures of the measure-valued branching diffusion with immigration  and the Fleming-Viot process with parent independent mutation. We view the relation between these two classes of measure-valued processes
      as the dynamical Gamma-Dirichlet structure. Other results of this article include the derivation of the transition function of the Fleming-Viot process with parent independent mutation from the transition function of the measure-valued branching diffusion with immigration, and the establishment of the reversibility of the latter. One of these is related to an open problem by Ethier and Griffiths and the other leads to an alternative proof of the reversibility of the Fleming-Viot process.

      \vspace*{.125in} \noindent {\it Keywords:} branching process with immigration,  coalescent, Dirichlet process, gamma process, Hamiltonian, large deviations,  quasi-invariant,  reversibility,
      random time-change.

      \vspace*{.125in} \noindent {\it AMS 2001 subject classifications:}
      Primary 60F10; secondary 92D10.
      \end{abstract}
      \section{Introduction}

      Recall that for any $\alpha >0, \beta >0$, the $Gamma(\alpha, \beta)$ distribution has the density function
\[
f(x)=\frac{1}{\Gamma(\alpha)\beta^{\alpha}}x^{\alpha-1}e^{-\frac{x}{\beta}}, \ x >0,
\]
and the Laplace transform
\[
\int_0^{\infty}e^{-u x}f(x)d\,x=\exp\{-\alpha\log(1+\beta u)\},\  u >-1/\beta,
\]
where $\alpha$ is the shape parameter and $\beta$ is the scale parameter. A characterization of the gamma distribution obtained in  \cite{Luc55} states that two independent positive random variables $Y_1$ and $Y_2$ are gamma random variables with the same scale parameter if and only if $Y_1+Y_2$ and $\frac{Y_1}{Y_1+Y_2}$ are independent.

For $i=1,2,$ let $Y_i$ be a $Gamma(\alpha_i, \beta)$ random variable with $\alpha_i>0$ and assume that $Y_1$ and $Y_2$ are independent. Then we have

\beqn
(Y_1,Y_2)&=& (Y_1+Y_2)(\frac{Y_1}{Y_1+Y_2}, \frac{Y_2}{Y_1+Y_2})\label{gb1}\\
&=&\mbox{radial part}\times \mbox{angular part},\nn
\eeqn
 where
 $Y_1+Y_2 $ and $(\frac{Y_1}{Y_1+Y_2}, \frac{Y_2}{Y_1+Y_2})$ are independent with respective distributions $Gamma(\alpha_1 +\alpha_2,\beta)$ and $Beta(\alpha_1,\alpha_2)$. The equation \rf{gb1} is called the {\it one-dimensional Gamma-Dirichlet structure}.

Let $S$ be a compact metric space, $\nu_0$ a probability on $S$, $\theta$ and $\beta$ any two positive numbers. The space of all non-negative finite measures on $S$, denoted by $M(S)$, is equipped with the weak topology, and $M_1(S)$, a subspace of $M(S)$, consists of all probability measures on $S$. The integration of a measurable function $g$ on $S$ with respect to a measure $\mu$ in $M(S)$ is denoted by $\langle \mu, g \rangle$. We denote by $B(S)$ and $C(S)$ the respective sets of bounded measurable functions and continuous functions on $S$. The {\it gamma process} with shape parameter $\theta \nu_0$ and scale parameter $\beta$ is given by

$$
{\cal Y}_{\theta,\nu_0}^{\beta}(\cdot)=\beta\sum_{i=1}^{\infty}\ga_i(\theta)\delta_{\xi_i}(\cdot)
$$
where $\ga_1(\theta)>\ga_2(\theta)> \cdots$ are the points of the inhomogeneous Poisson point process on $(0,\infty)$ with mean measure $\theta x^{-1}e^{-x}d\,x$, and independently, $\xi_1,\xi_2,\ldots$ are i.i.d. with common distribution $\nu_0$.

Denote the law of ${\cal Y}_{\theta,\nu_0}^{\beta}$ by $\Gamma_{\theta,\nu_0}^{\beta}$. The corresponding Laplace functional is

\[
\int_{M(S)}e^{-\langle \mu, g \rangle}\Gamma_{\theta,\nu_0}^{\beta}(d\,\mu) =\exp\{-\theta \langle \nu_0, \log(1+\beta g)\rangle\}
\]
where
\[
g(s)>-1/\beta, \ \mb{for all }\ s \in S.
\]

Set
\[
 \sigma(\theta)=\sum_{i=1}^{\infty}\ga_i(\theta),\ \ P_i(\theta) =\frac{\ga_i(\theta)}{\sigma(\theta)}, \ \ {\cal X}_{\theta, \nu_0}(\cdot)=\sum_{i=1}^{\infty}P_i(\theta) \delta_{\xi_i}(\cdot).
\]
\vspace{0.3cm}

The law of $(P_1(\theta),P_2(\theta),\ldots)$, denoted by $PD(\theta)$, is the {\it Poisson-Dirichlet distribution} with parameter $\theta$ (cf. \cite{Kingman75}), ${\cal X}_{\theta, \nu_0}(\cdot)$ is the {\it Dirichlet process} with law denoted by $\Pi_{\theta,\nu_0}$, and the relation
\be\label{gb2}
{\cal X}_{\theta, \nu_0}(\cdot)= \frac{{\cal Y}_{\theta,\nu_0}^{\beta}(\cdot)}{{\cal Y}_{\theta,\nu_0}^{\beta}(S)}
\ee
is called the {\it infinite-dimensional Gamma-Dirichlet structure}.

The law $\Gamma^{\beta}_{\theta,\nu_0}$ is the reversible measure of a measure-valued branching diffusion with immigration (henceforth MBI) with generator
\be\label{gen-mbi}
{\cal L}=\frac{1}{2}\{\langle \mu, \frac{\delta^2}{\delta \mu(s)^2}\rangle + \langle \theta\nu_0-\la \mu, \frac{\delta}{\delta \mu(s)} \rangle\}
\ee
where $\la =\frac{1}{\beta}>0$ and
\[
\frac{\delta \varphi}{\delta \mu(s)}=\lim_{\ep \ra 0}\frac{\varphi(\mu +\ep \delta_s)-\varphi(\mu)}{\ep}.
\]
The law $\Pi_{\theta,\nu_0}$ is the reversible measure of the Fleming-Viot process with parent independent mutation (henceforth FVP) with generator
\be\label{gen-fvp}
{\cal A}=\frac{1}{2}\{\langle \nu(d\,s_1)\cdot(\delta_{s_1} -\nu)(d\,s_2), \frac{\delta^2}{\delta \nu(s_2)\delta\nu(s_1)}\rangle + \theta\langle \nu_0- \nu, \frac{\delta}{\delta \nu(s)} \rangle\}.
\ee
The {\it dynamical Gamma-Dirichlet structure} corresponds to the relation between these two classes of measure-valued processes.

Section 2 reviews some known results including several algebraic identities, the formal Hamiltonian, the quasi-invariant, and large deviations. The large deviation principle (henceforth LDP) for the jump sizes of the gamma process is established in Section 3. Section 4 obtains new results on quasi-invariance for the jump size of gamma process and the two-parameter Poisson-Dirichlet distribution. Finally in section 5, we derive the transition function of the FVP process directly from the transition function of the MBI process through a time change, and establish the reversibility of the MBI process. This, combined with the Gamma-Dirichlet structure, provides an alternative proof of the reversibility of the FVP process.

      \section{Gamma-Dirichlet Structure}
      \setcounter{equation}{0}

   The gamma distribution is structurally similar to the normal distribution in many ways. For example, a normal population is characterized by the independency of any translation-invariant statistic and the sample mean while a gamma population is characterized by the independency of any scale-invariant statistic of the sample mean or equivalently by the one-dimensional Gamma-Dirichlet structure. In this section, we collect several existing results associated with the Gamma-Dirichlet structure, which motivate the studies in subsequent sections.

\subsection{Algebraic Relations}
 Let $\overset{d}{=}$ denote equality in distribution. For any integer $m \geq 2$, $\beta >0$, $\theta_i>0$, $\nu_i \in M_1(S),\linebreak i=1,\ldots,m$, let ${\cal Y}^{\beta}_{\theta_1,\nu_1}, \ldots, {\cal Y}_{\theta_m, \nu_m}^{\beta}$ be $m$ independent gamma processes with the same scale parameter $\beta$.  By direct calculation, we obtain the following additive property:

\be \label{add-pro1}
{\cal Y}^{\beta}_{\theta_1,\nu_1} +\cdots+{\cal Y}_{\theta_m, \nu_m}^{\beta} \overset{d}{=} {\cal Y}^{\beta}_{\theta_1+\cdots+\theta_m, \frac{\theta_1}{\theta_1+\cdots+\theta_m}\nu_1+\cdots+\frac{\theta_m}{\theta_1+\cdots+\theta_m}\nu_m}.
\ee
\vspace{0.1cm}

Let ${\cal X}_{\theta_1, \nu_1}, \ldots, {\cal X}_{\theta_m,\nu_m}$ be $m$ independent Dirichlet processes, and independently let the random vector $(\eta_1,\ldots,\eta_m)$ have a $Dirichlet(\theta_1,\ldots, \theta_m)$ distribution. Then the following mixing identity follows from \rf{add-pro1} and the infinite-dimensional Gamma-Dirichlet structure.

\be\label{add-pro2}
\eta_1 {\cal X}_{\theta_1,\nu_1}+\cdots +\eta_m{\cal X}_{\theta_m, \nu_m} \overset{d}{=} {\cal X}_{\theta_1+\cdots+\theta_m, \frac{\theta_1}{\theta_1+\cdots+\theta_m}\nu_1+\cdots+\frac{\theta_m}{\theta_1+\cdots+\theta_m}\nu_m}.
\ee

Let $U_1(\theta), U_2(\theta), \ldots$ be a sequence of i.i.d. random variables with common distribution $Beta(1,\theta)$. Set
\be\label{gem}
V_1(\theta)=U_1(\theta),\ V_n(\theta) = U_n(\theta)\prod_{i=1}^{n-1}(1-U_i(\theta)),\ n \geq 2.
\ee

Then the law of $(V_1(\theta), V_2(\theta),\ldots)$ is the GEM distribution with parameter $\theta$ and the law of the descending order statistics of $(V_1(\theta),V_2(\theta),\ldots)$ is the Poisson-Dirichlet distribution $PD(\theta)$. The Dirichlet process ${\cal X}_{\theta, \nu_0}$ is then given by
\[
{\cal X}_{\theta, \nu_0}=\sum_{i=1}^{\infty}V_i(\theta) \delta_{\xi_i},
\]
where $(V_1(\theta),V_2(\theta),\ldots)$ is independent of $(\xi_1,\xi_2,\ldots)$. It follows from the  self-similarity of the GEM representation \rf{gem} that

\be\label{diri1}
{\cal X}_{\theta, \nu_0}= V_1(\theta) \delta_{\xi_1}+(1-V_1(\theta))\tilde{{\cal X}}_{\theta, \nu_0},
\ee
where $\tilde{{\cal X}}_{\theta, \nu_0}$ is an independent copy of ${\cal X}_{\theta, \nu_0}$.

\subsection{Quasi-Invariance}

 Consider a measure $P$ on space $S$. Let $G$ denote a group of transformations from $S$ to $S$. The measure $P$ is {\it quasi-invariant}
with respect to the group $G$ if for any $T$ in $G$ the image measure $P\circ T^{-1}$ of $P$ under $T$ is equivalent to $P$. In this case, $T$ preserves zero sets. If $P\circ T^{-1}=P$ , the measure $P$ is {\it invariant} with respect to $G$.

The quasi-invariant property of Gaussian measure plays a major role in stochastic calculus. In \cite{TV99} and \cite{TVY01}, the quasi-invariance of the gamma process is studied thoroughly. A comparison with the Gaussian measure reveals a remarkable similarity between the spherical symmetry of the Gaussian measure and the multiplicative symmetry of the gamma process. This provides a different aspect for the Gamma-Dirichlet structure.

Let $B_+(S)$ be the collection of positive Borel measurable functions on $S$ with strictly positive lower bound. For each $f$ in $B_+(S)$, set
\[
T_f(\mu)(d\,s)=f(s)\mu(d\,s),\  \mu \in M(S),
\]
and let $T_f(\Gamma_{\theta,\nu_0}^{\beta})$ denote the image law of $\Gamma_{\theta,\nu_0}^{\beta}$ under $T_f$. Then the following holds.

\begin{theorem}\label{t1}{\rm (Tsilevich and Vershik \cite{TV99}, Tsilevich, Vershik and Yor \cite{TVY01})}.
For any $f$ in $B_+(S)$,
$T_f(\Gamma_{\theta,\nu_0}^{\beta})$ and $\Gamma_{\theta,\nu_0}^{\beta}$ are mutually absolutely continuous and

\be\label{qu-den1}
\frac{d\, T_f(\Gamma_{\theta,\nu_0}^{\beta})}{d\, \Gamma_{\theta,\nu_0}^{\beta}}(\mu)=\exp\{-[\theta\langle \nu_0, \log f\rangle +\langle \mu, \beta^{-1}(f^{-1}-1)]\}.
\ee
\end{theorem}

\vspace{0.3cm}

For each $f$ in $B_+(S)$, set
\[
\overline{T}_f(\nu)(d\,s)=\frac{f(s)\nu(d\,s)}{\langle \nu, f\rangle },\  \nu \in M_1(S),
\]
and denote by $\overline{T}_f(\Pi_{\theta,\nu_0})$ the image law of $\Pi_{\theta,\nu_0}$ under $\overline{T}_f$.

\begin{theorem}{\rm (Handa \cite{Handa01})}\label{diri2}
The laws $\overline{T}_f(\Pi_{\theta,\nu_0})$ and $\Pi_{\theta,\nu_0}$ are mutually absolutely continuous and

\be\label{qu-den2}
\frac{d\, \overline{T}_f(\Pi_{\theta,\nu_0})}{d\, \Pi_{\theta,\nu_0}}(\nu)=\exp\{-\theta[\langle \nu_0, \log f\rangle +\log\langle \nu, f^{-1}\rangle]\}.
\ee
\end{theorem}

\vspace{0.3cm}

Let
\[
\nabla_{\infty}=\left\{{\bf p}=(p_1,p_2,\ldots): p_1\geq p_2\geq\cdots\geq 0, \sum_{j=1}^{\infty}p_j=1\right\}
\]
be equipped with the subspace topology of the infinite-dimensional Euclidean space. Then the Poisson-Dirichlet distribution $PD(\theta)$ is a probability measure on $\nabla_{\infty}$. For any $f$ in $B_+(S)$ and any i.i.d. sequence $\{\xi_i: i=1,2,\ldots\}$ with common distribution $\nu_0$, define the map
\be\label{atom-map}
S_{f,\nu_0}: \nabla_{\infty} \ra \nabla_{\infty}, \ (p_1,p_2,\ldots)\mapsto (\tilde{p}_1,\tilde{p}_2, \ldots)
\ee
where $(\tilde{p}_1,\tilde{p}_2, \ldots)$ is the descending order statistics of
\[
\left(\frac{f(\xi_1)p_1}{\sum_{i=1}^{\infty}f(\xi_i)p_i}, \frac{f(\xi_2)p_2}{\sum_{i=1}^{\infty}f(\xi_i)p_i},\ldots\right).
\]

\begin{theorem}\label{diri4}{\rm ( Tsilevich and Vershik \cite{TV99})}
Let $\widetilde{PD}(\theta)$ denote the image law of $PD(\theta)$ under the map $S_{f,\nu_0}$. Then $\widetilde{PD}(\theta)$ and $PD(\theta)$ are mutually absolutely continuous and

\be\label{qu-den3}
\frac{d\, \widetilde{PD}(\theta)}{d\, PD(\theta)}({\bf p})=\exp\{-\theta\langle \nu_0, \log f\rangle\}\cdot \int_0^{\infty}\frac{u^{\theta-1}}{\Gamma(\theta)}\prod_{i=1}^{\infty}\langle \nu_0, e^{-u f^{-1}(\xi_i)p_i}\rangle d\,u.
\ee
\end{theorem}

\vspace{0.3cm}

Consider an abstract space $\Omega$ with a formal reference probability measure ${\cal P}$ (uniform or invariant in some sense). The {\it formal Hamiltonian} ${\cal H}(\omega)$ is a function associated with another probability ${\cal Q}$ such that
\[
{\cal Q}(d\omega)=Z^{-1}\exp\{-{\cal H(\omega)}\}{\cal P}(d\,\omega).
\]

For each $\mu \in M(S)\setminus \{0\}$, set $\hat{\mu}(\cdot)=\frac{\mu(\cdot)}{\mu(S)} \in M_1(S)$. Let
\[
\phi(x)=x\log x -(x -1), x \geq 0,
\]
and for any $\nu_1, \nu_2 \in M_1(S)$
\[
Ent(\nu_1|\nu_2)=\left\{\begin{array}{ll}
     \int_{M_1(S)}\log\frac{d\,\nu_1}{d\,\nu_2}d\,\nu_1 ,  & \nu_1\ll\nu_2\\
     +\infty,  & \mbox{else}.
      \end{array}\right.
\]

Then the following holds.

\begin{theorem}\label{diri3}{\rm (Handa \cite{Handa01})}
The formal Hamiltonian for the Gamma process $\Gamma_{\theta,\nu_0}^{\beta}$ is given by
\beqn\label{g-qua}
{\cal H}_g(\mu) &=&\theta Ent(\nu_0|\hat{\mu})+ \frac{\mu(S)}{\beta}\phi(\frac{\beta\theta}{\mu(S)})\\
&=& \mb{angular component} +\mb{radial component},\nn
\eeqn
and the formal Hamiltonian for the Dirichlet process $\Pi_{\theta,\nu_0}$ is given by
\be\label{d-qua}
{\cal H}_d(\nu) =\theta Ent(\nu_0|\nu), \ \nu \in M_1(S).
\ee
\end{theorem}

\vspace{0.3cm}

\noindent {\bf Remark} For $\beta \theta=1$, we have
\[
{\cal H}_g(\mu)|_{\mu(S)=1} ={\cal H}_d(\hat{\mu})
\]
which, combined with \rf{g-qua}, provides a different version of the Gamma-Dirichlet structure.

\section{ Asymptotic Behavior for the Jump Sizes of the Gamma Process}
\setcounter{equation}{0}

The Gamma-Dirichlet structure has been used in \cite{Gri79} and \cite{Feng10} to obtain the fluctuation theorem for the Poisson-Dirichlet distribution. In this section, we apply the Gamma-Dirichlet structure to the establishment of the LDP for the jump sizes of the Gamma process.

Large deviations for a family of probability measures $\{P_{\la}: \la \in \mb{index
 set}\}$ on a space $E$ are estimations of
 the following type:
 \[
   P_{\la}\{G\} \asymp \exp\{-a(\la)\inf_{x \in G}I(x) \},
 \]
 where $a(\la)$ is called the large deviation speed, and the nonnegative lower semi-continuous function $I(\cdot)$ is called the rate function. The rate function is good if $\{x\in E: I(x)\leq c \}$ is compact for any nonnegative constant $c$. Additional terminologies and general results on LDPs are found in \cite{DZ98}.

 Recall that $\ga_1(\theta)>\ga_2(\theta)> \cdots$ are the points of the inhomogeneous Poisson point process on $(0,\infty)$ with mean measure $\theta x^{-1}e^{-x}d\,x$, and

\be\label{gd-atom}
(\ga_1(\theta),\ga_2(\theta),\ldots)=\sigma(\theta) (P_1(\theta),P_2(\theta),\ldots)
\ee
where $\sigma(\theta)$ is a $Gamma(\theta,1)$ random variable that is independent of the Poisson-Dirichlet distributed random sequence $(P_1(\theta),P_2(\theta),\ldots)$.

\begin{theorem}\label{ldp-gama}
Let $P_{\theta}$ denote the law of $\theta^{-1}(\ga_1(\theta), \ga_2(\theta), \ldots)$ on $\mathbb{R}_+^{\infty}$. Then the family $\{P_{\theta}: \theta >0\}$ satisfies a LDP with speed $\theta$ and good rate function
\[
I(x_1,x_2,\ldots)=\left\{\begin{array}{ll}
\displaystyle \sum_{i=1}^{\infty}x_i,& x_1\geq x_2\geq\cdots\geq 0 \\
\displaystyle\infty,& \mb{otherwise},
\end{array}\right.
\]
as $\theta$ converges to infinity.
\end{theorem}
\proof It follows from direct calculation that the laws of the family $\{\theta^{-1}\sigma(\theta): \theta >0\}$ satisfy a LDP as $\theta$ tends to infinity with speed $\theta$ and good rate function

\[
  I_1(y)=\left\{\begin{array}{ll}
\displaystyle y-1-\log y,& y>0 \\
\displaystyle\infty,& \mb{otherwise},
\end{array}\right.
 \]

Let

\[
      \nabla = \{(p_1,...,p_n,...): p_1\geq p_2\geq\cdots\geq 0,\  \sum_{i=1}^{\infty} p_i \leq 1\}.
\]

It follows from Theorem 4.4 in \cite{DF05} that the family $\{PD(\theta): \theta >0\}$ satisfies a LDP on space $\nabla$ as $\theta $ tends to infinity with speed $\theta$ and good rate function

\[
 I_2(p_1,p_2,\ldots)=\left\{\begin{array}{ll}
\displaystyle -\log (1-\sum_{i=1}^{\infty}p_i),& \sum_{i=1}^{\infty}p_i <1 \\
\displaystyle\infty,& \mb{otherwise},
\end{array}\right.
 \]

The Gamma-Dirichlet structure \rf{gd-atom}, combined with the contraction principle, implies that the family $\{P_{\theta}: \theta >0\}$ satisfies a LDP as $\theta $ tends to infinity with speed $\theta$ and good rate function
 \[
 \inf\{I_1(y)+I_2(p_1,p_2,\ldots): y\geq 0, (p_1,\ldots)\in \nabla,  x_i=y p_i, i=1,\ldots\}
 \]
 which equals to $I(x_1,x_2,\ldots)$.

 \hfill $\Box$

 \noindent {\bf Remarks}  The LDP  is established in \cite{Feng09} for $PD(\theta)$ when $\theta$ converges to zero. The speed is $-\log(\theta)$ and the good rate function is
 \[
 I_3(p_1,p_2,\ldots)=\left\{\begin{array}{ll}
      0,& p_1 =1\\
      n-1,&  \sum_{i=1}^{n}p_i =1, p_n >0\\
      \infty,& \mb{else}.\end{array}\right.
 \]

 Since $(\ga_1(\theta), \ga_2(\theta), \ldots)$ converges to $(0,0,\ldots)$ when $\theta$ converges to zero, one would like to establish the LDP for $(\ga_1(\theta), \ga_2(\theta), \ldots)$ from the LDP for $PD(\theta)$. A direct calculation shows that a LDP  holds for $\sigma(\theta)$ with speed $-\log(\theta)$ and rate function
 \[
  I_4(y)=\left\{\begin{array}{ll}
\displaystyle 0,& y=0 \\
\displaystyle 1,& \mb{otherwise},
\end{array}\right.
 \]
Since  $I_4(\cdot)$ is not a good rate function, the contraction principle can not be applied in this case. On the other hand, for any $r \geq 1$ the joint density function of $(\ga_1(\theta), \ldots, \ga_r(\theta))$ is
\[
\frac{\theta^r}{x_1\cdots x_r}\exp\{-\sum_{i=1}^r x_i-\theta E_1(x_r)\}, \ \  x_1\geq \cdots\geq x_r >0
\]
where
\[
E_1(x)=\int_x^{\infty}z^{-1}e^{-z}d\,z.
\]
By direct calculation, one can show that a LDP holds for  $(\ga_1(\theta),\ga_2(\theta),\ldots)$ with speed $-\log(\theta)$ and rate function
\[
 I_5(x_1,x_2,\ldots)=\left\{\begin{array}{ll}
\displaystyle 0,& x_i=0, i\geq 1 \\
\displaystyle n,& x_1\geq \cdots \geq x_n >0, x_m=0, m\geq n+1 \\
\displaystyle \infty,& \mb{otherwise},
\end{array}\right.
\]
Clearly $I_5(\cdot)$ is not a good rate function. But remarkably the three rate functions have a relation that is consistent with the contraction principle, i.e.,
\[
I_5(x_1,x_2,\ldots)=\inf\{I_3(p_1,\ldots) +I_4(y): yp_i=x_i, i=1,\ldots\}.
\]

 \section{New Results on Quasi-Invariance}
\setcounter{equation}{0}

For any $0<\alpha<1$ and $\theta >-\alpha$, let $U_k(\alpha,\theta),
      k=1,2,...$, be a sequence of independent
      random variables such that $U_i(\alpha,\theta)$ has $Beta(1-\alpha,\theta+ i\alpha)$ distribution.
      Set
      \be \label{GEM2}
      V_1(\alpha,\theta) = U_1(\alpha,\theta),\  V_n(\alpha,\theta) = (1-U_1(\alpha,\theta))\cdots (1-U_{n-1}(\alpha,\theta))U_n(\alpha,\theta),\  n \geq 2.
      \ee
      Then the law of $(V_1(\alpha,\theta),V_2(\alpha, \theta),...)$
      is called the two-parameter GEM distribution, denoted by $GEM(\theta,\alpha).$ Let ${\bf P}(\alpha,\theta)=(P_1(\alpha,\theta), P_2(\alpha,\theta),...)$ denote $(V_1(\alpha,\theta),V_2(\alpha,\theta),...)$
      in descending order. The law of ${\bf P}(\alpha,\theta)$ is called the two-parameter Poisson-Dirichlet
      distribution, and, following \cite{PitmanYor97}, is denoted by
      $PD(\alpha,\theta)$. In this section, we will generalize Theorem~\ref{t1} and Theorem~\ref{diri4} to the jump sizes of the gamma process and the two-parameter Poisson-Dirichlet distribution, respectively.

  Let

  \[
  M_a(S)= \{\sum_{i=1}^{\infty}x_i\delta_{\xi_i}: \xi_i \in S, x_i \geq 0, 0<\sum_{i=1}^{\infty}x_i < \infty\}
  \]
denote the space of non-negative finite atomic measures on $S$, equipped with the subspace topology of $M(S)$.
 Let
 \[
 \mathbb{R}_+^\downarrow=\{{\bf x}=(x_1,x_2,\ldots): x_1\geq x_2\geq\cdots\geq 0, \sum_{i=1}^{\infty}x_i <\infty\},
 \]
 and for any $\nu=\sum_{i=1}^{\infty}y_i\delta_{\xi_i}$ in
$M_a(S)$ define the map
\[
{\cal J}: M_a(S)\ra \mathbb{R}_+^\downarrow, \nu\mapsto {\bf x},
\]
where ${\bf x} $ is $(y_1,y_2,\ldots)$ in descending order. Similarly we define
 the map
\[
{\cal T}: M_a(S)\ra \nabla_{\infty},\ \nu \mapsto {\bf p}=(p_1,p_2,\ldots)
\]
where ${\bf p}$ is $\{\frac{y_i}{\sum_{j=1}^{\infty}y_j}: i=1,2,\ldots\}$ in descending order. Since the usual weak topology generates the same Borel $\sigma$-field on $M_a(S)$ as the weak atomic topology, it follows from Lemma 2.5 in \cite{EtKu94} that the maps ${\cal J}$ and ${\cal T}$ are measurable.

\begin{theorem}\label{atom-gamma}
For any $f$ in $B_+(S)$, let $\tilde{T}_f(\Gamma_{\theta,\nu_0}^{\beta})$ and $\tilde{\Gamma}_{\theta,\nu_0}^{\beta}$ be the respective  image laws of
$T_f(\Gamma_{\theta,\nu_0}^{\beta})$ and $\Gamma_{\theta,\nu_0}^{\beta}$ under ${\cal J}$. Then $\tilde{T}_f(\Gamma_{\theta,\nu_0}^{\beta})$ and $\tilde{\Gamma}_{\theta,\nu_0}^{\beta}$ are mutually absolutely continuous and

\be\label{atom-gamma-ac}
\frac{d\, \tilde{T}_f(\Gamma_{\theta,\nu_0}^{\beta})}{d\, \tilde{\Gamma}_{\theta,\nu_0}^{\beta}}({\bf x})=
\exp\{-\theta\langle \nu_0, \log f\rangle\}E^{\nu_0^{\infty}}
[\exp\{-\sum_{i=1}^{\infty}\beta^{-1}(f^{-1}(\xi_i)-1)x_i\}].
\ee
\end{theorem}

\proof For any bounded measurable function $F$ on $\mathbb{R}_+^{\downarrow}$, it follows from Theorem~\ref{t1} that
\beq
&&E^{\tilde{T}_f(\Gamma_{\theta,\nu_0}^{\beta})}[F({\bf x})]= E^{T_f(\Gamma_{\theta,\nu_0}^{\beta})}[F({\cal J}(\mu))]\\
&&\hspace{1cm}= \exp\{-\theta\langle \nu_0, \log f\rangle\}E^{\Gamma_{\theta,\nu_0}^{\beta}}[F({\cal J}(\mu))\exp\{-\langle \mu, \beta^{-1}(f^{-1}-1)\}]\\
&&\hspace{1cm}= \exp\{-\theta\langle \nu_0, \log f\rangle\}E^{\tilde{\Gamma}_{\theta,\nu_0}^{\beta}}\left[F({\bf x})E^{\nu_0^{\infty}}[\exp\{-\sum_{i=1}^{\infty}\beta^{-1}(f^{-1}(\xi_i)-1)x_i\}]\right]
\eeq
which leads to the result.

\hfill $\Box$

The stable subordinator with index $\alpha$ is a L\'{e}vy process with L\'{e}vy measure
$$
d\Lambda_\alpha=\frac{c\alpha}{\Gamma(1-\alpha)}s^{-\alpha-1}d\,s, \  s>0
$$
for some $c>0$. Let $\rho_1(\alpha)\geq \rho_2(\alpha)\geq \ldots$ denote the descending order jump sizes of the stable subordinator over the interval $(0,1]$. Independently, let $\xi_1, \xi_2, \ldots$ be i.i.d. with common distribution $\nu_0$. Then the $\nu_0$-scaled stable subordinator has the form
\[
\sum_{i=1}^{\infty}\rho_i(\alpha)\delta_{\xi_i}
\]
with law denoted by $P_{\alpha,0}^{c,\nu_0}$.
The Laplace functional of $P_{\alpha,0}^{c,\nu_0}$
is given by
$$
E^{P_{\alpha,0}^{c,\nu_0}}[\exp\{-\langle \nu, g\rangle\}]=\exp\{-c\langle \nu_0,g^\alpha\rangle\}.
$$

\begin{lemma}\label{new-qua}
For any $f$ in $B_+(S)$, let $A_\alpha=\langle \nu_0,f^\alpha\rangle$ and $\nu_\alpha(d\,s)=f^\alpha(s)\nu_0(d\,s)/A_\alpha.$ Then the law $T_f(P_{\alpha,0}^{c,\nu_0})$
is $P_{\alpha,0}^{c A_\alpha,\nu_\alpha}$.
\end{lemma}
\proof For any $g$ in $B_+(S)$, it follows from direct calculation that
\begin{eqnarray*}
E^{T_f(P_{\alpha,0}^{c,\nu_0})}[\exp\{-\langle \nu, g \rangle\}]&=&E^{P_{\alpha,0}^{c,\nu_0}}[\exp\{-\langle \nu,gf\rangle\}]\\
&=&\exp\{-c\langle \nu_0, f^\alpha g^\alpha\rangle\}\\
&=&\exp\left\{-c A_\alpha\langle \frac{f^\alpha \nu_0}{A_\alpha}, g^\alpha, \rangle\right\}\\
&=&E^{P_{\alpha,0}^{cA_\alpha,\nu^\alpha_0}}[\exp\{-\langle \nu, g\rangle\}],
\end{eqnarray*}
which leads to the result.

\hfill $\Box$

Next set $c_{\alpha,\theta}=c^{\theta/\alpha}\frac{\Gamma(\theta+1)}{\Gamma(\theta/\alpha+1)}$ and consider the law $P^{c,\nu_0}_{\alpha,\theta}$ defined by
\be\label{gir1}
P^{c,\nu_0}_{\alpha,\theta}(d\,\nu)=\frac{c_{\alpha,\theta}}{\nu(S)^\theta}P^{c,\nu_0}_{\alpha,0}(d\,\nu).
\ee

It is known (cf. \cite{PPY92}, \cite{PitmanYor97}) that
\be\label{map1}
{\cal T}(P^{c,\nu_0}_{\alpha,0})=PD(\alpha,0)\ee
 and
 \be\label{map2}
 {\cal T}(P^{c,\nu_0}_{\alpha,\theta})=PD(\alpha,\theta).
 \ee

The next theorem generalizes the result in Theorem~\ref{diri4} to the two-parameter Poisson-Dirichlet distribution $PD(\alpha, \theta)$.

\begin{theorem}\label{diri5}
For any $f$ in $B_+(S)$, let $\widetilde{PD}(\alpha,\theta)$ denote the image law of $PD(\alpha,\theta)$ under the map $S_{f,\nu_0}$ defined in \rf{atom-map}. Then $\widetilde{PD}(\alpha, \theta)$ and $PD(\alpha,\theta)$ are mutually absolutely continuous and

\be\label{qu-den4}
\frac{d\,\widetilde{PD}(\alpha,\theta)}{d\,PD(\alpha,\theta)}({\bf p})=\langle \nu_0,f^\alpha\rangle^{-\theta/\alpha}\frac{1}{\Gamma(\theta)}\int_0^\infty u^{\theta-1}\prod_{i=1}^\infty E^{\nu_\alpha}[e^{-u f^{-1}(\xi_i)p_i}] d\,u.
\ee
\end{theorem}
\proof For any nonnegative product measurable function $\Phi$ on $\nabla_{\infty}$, we have
\beq
E^{\widetilde{PD}(\alpha,\theta)}[\Phi({\bf p})]&=&E^{T_f(P^{c,\nu_0}_{\alpha,\theta})}[\Phi({\cal T}(\mu))]\\
&=&E^{P^{c,\nu_0}_{\alpha,\theta}}[\Phi({\cal T}(T_f(\nu))]\\
&=&E^{P_{\alpha,0}^{c,\nu_0}}\left[\Phi({\cal T}(T_f(\nu))\frac{c_{\alpha,\theta}}{\nu(S)^\theta}\right].
\eeq

Since $\nu(S)$ can be written as $\langle T_f(\nu), f^{-1}\rangle$, it follows from lemma~\ref{new-qua} that
\beq
E^{P_{\alpha,0}^{c,\nu_0}}\left[\Phi({\cal T}(T_f(\nu))\frac{c_{\alpha,\theta}}{\nu(S)^\theta}\right]
&=&E^{P_{\alpha,0}^{c,\nu_0}}\left[\Phi({\cal T}(T_f(\nu))\frac{c_{\alpha,\theta}}{\langle T_f(\nu), f^{-1}\rangle^\theta}\right]\\
&=&E^{T_f(P_{\alpha,0}^{c,\nu_0})}\left[\Phi({\cal T}(\mu))\frac{c_{\alpha,\theta}}{\langle \mu,f^{-1}\rangle^\theta}\right]\\
&=&E^{T_f(P_{\alpha,0}^{c A_{\alpha},\nu_{\alpha}})}\left[\Phi({\cal T}(\mu))\frac{c_{\alpha,\theta}}{\langle \mu,f^{-1}\rangle^\theta}\right]\\
&=&E^{P_{\alpha,\theta}^{cA_\alpha,\nu_\alpha}}\left[\Phi({\cal T}(\mu))\frac{c_{\alpha,\theta}}{(cA_\alpha)_{\alpha,\theta}}
\frac{\mu(S)^\theta}{\langle \mu,f^{-1}\rangle^\theta}\right]\\
&=&E^{PD(\alpha,\theta)}\left[\Phi({\bf p})E^{\nu_\alpha^\infty}[\frac{1}{A_\alpha^{\theta/\alpha}(\sum f^{-1}(\xi_i)p_i)^\theta}]\right],
\eeq
where $\nu_\alpha^\infty$ is the infinite product of $\nu_{\alpha}$.
Hence we have
$$
\frac{\widetilde{PD}(\alpha,\theta)}{PD(\alpha,\theta)}(d{\bf p})=E^{\nu_\alpha^\infty}[\frac{1}{A_\alpha^{\theta/\alpha}(\sum f^{-1}(\xi_i)p_i)^\theta}],
$$

Note that for $\lambda > 0$, $\lambda^{-\theta}=\int_0^\infty \frac{\sigma^{\theta-1}}{\Gamma(\theta)}e^{-\lambda\sigma}d\sigma$. Consequently,
\begin{eqnarray*}
E^{\nu_\alpha^\infty}\left[\frac{1}{(\sum f^{-1}(\xi_i)p_i)^\theta}\right]&=&E^{\nu_\alpha^\infty}\left[\int_0^\infty
\frac{\sigma^{\theta-1}}{\Gamma(\theta)}e^{-\sigma\sum f^{-1}(\xi_i)p_i}d\sigma\right]\\
&=&\frac{1}{\Gamma(\theta)}\int_0^\infty \sigma^{\theta-1}\prod_{i=1}^\infty E^{\nu_\alpha}[e^{-\sigma f^{-1}(\xi_i)p_i}] d\sigma,
\end{eqnarray*}
which leads to \rf{qu-den4}. The fact that $f$ has a strictly positive lower bound implies that $$\prod_{i=1}^\infty E^{\nu_\alpha}[e^{-\sigma f^{-1}(\xi_i)p_i}]$$ is strictly positive. Therefore the right hand side of \rf{qu-den4} is
finite everywhere.

\hfill $\Box$

\noindent {\bf Remarks} When $\theta=0$, the density equals to one and $P_{\alpha,0}$ is a fixed point for the map $S_{f,\nu_0}$. On the other hand, by direct calculation it follows that
\beq
\nu_{\alpha}&\Rightarrow&\nu_0, \ \alpha \ra 0\\
\lim_{\alpha\ra 0}\langle \nu_0, f^\alpha\rangle^{-\theta/\alpha}&=&\exp\{-\theta\langle \nu_0, \log f\rangle\}.
\eeq
Therefore, for positive $\theta$, the quasi-invariance of $PD(\theta)$ can be obtained from the quasi-invariance of $PD(\alpha,\theta)$ by taking the limit of $\alpha$ going to zero.

\section{Dynamical Gamma-Dirichlet Structure}
\setcounter{equation}{0}

Let $Y_t$ denote the MBI process with generator ${\cal L}$ given by \rf{gen-mbi}. The process $X_t$ denotes the FVP process with generator ${\cal A}$ given by \rf{gen-fvp}.  In this section, we will explore the dynamical Gamma-Dirichlet structure between $Y_t$ and $X_t$.

\subsection{Transition Functions and Random Time Change}

For $\la \geq 0, t \geq 0$, set

\begin{equation}\label{CC}
C(\lambda,t)=\left\{\begin{array}{lr}\lambda^{-1}(\re^{\lambda
t/2}-1),&\lambda\not=0,\\t/2,&\lambda=0.\end{array}\right.
\end{equation}
\vspace{0.3cm}

For any $a>0$, let $\{N_a(t): t\geq 0\}$ be a time inhomogeneous pure death Markov
 chain with death rate $N_a(t)/2C(-\la,t)$ at time $t>0$, entrance boundary $\infty$, and marginal distribution

\[
q^{a,\la}_n(t)=P\{N_a(t)=n\}= \frac{a^n}{C^n(\la,t) n!}\exp\{-a/C(\la,t)\}, \ n =0,1,\ldots.
\]
\vspace{0.3cm}

 Let $Z_+ =\{0,1,\ldots\}$ and $\hat{Z}_+$ denote the one-point compactification of $Z_+$ by $\infty$. For any $\theta >0$, let $\{D_t^{\theta}: t \geq 0\}$ be the embedded chain of Kingman's
coalescent. It is a pure death process with state space $\hat{Z}_+$, death rates
\[
\la_n = \frac{n(n+\theta-1)}{2},\ \ n =0,1,\ldots,
\]
and entrance boundary $\infty$. For $t>0$, let
$$
d_n^\theta(t)=P(D_t^\theta=n)
$$
denote the probability of having $n$ lines of decent at time $t$  beginning at generation zero.

The following explicit formula for $d_n^\theta(t)$
is obtained in \cite{Tavar84}.
$$
d_n^\theta(t)=\left\{\begin{array}{lr} 1-\sum_{m=1}\limits^\infty
(2m-1+\theta)(m!)^{-1}(-1)^{m-1}\theta_{(m-1)}\re^{-\lambda_mt},&
n=0\\
\sum\limits_{m=n}^\infty(2m-1+\theta)(m!)^{-1}(-1)^{m-n}{m \choose
n}(n+\theta)_{(m-1)}\re^{-\lambda_m t},&  n\geq 1.
\end{array}\right.
$$

\begin{theorem}{\rm (Ethier and Griffiths \cite{EtGri93a}, \cite{EtGri93b})}\label{t1-sec5}
Assume that $t >0$ and  $\mu$ is in $M(S)$ with $\mu(S)>0$.  Set $\nu(\cdot)=\frac{\mu(\cdot)}{\mu(S)}.$

{\rm (1)} The transition function of the MBI process $Y_t$ is
\beqn
&&Q_{1}(t, \mu, \cdot)=q_0^{\mu(S),\la}\Gamma^{C(-\la,t)}_{\theta,\nu_0}(\cdot)\label{tran1}\\
&& \ \ \ \ \ \ +\sum_{n=1}^{\infty}q^{\mu(S),\la}_n(t)\int_{S^n}(\frac{\mu}{\mu(S)})^n(d\,x_1\times \cdots\times d\,x_n)
\Gamma^{C(-\la,t)}_{n+\theta, \frac{n}{\theta+n} \eta_n +\frac{\theta}{\theta+n}\nu_0}(\cdot)\nn
\eeqn
where $\eta_n = \frac{1}{n}\sum_{i=1}^n \delta_{x_i}$. The probability $\Gamma_{\theta,\nu_0}^{\beta}$ is the reversible measure for $Y_t$.

{\rm (2)} The transition function of the FVP process $X_t$ is
\beqn
&&Q_{2}(t, \nu, \cdot)=d_0^{\theta}(t) \Pi_{\theta,\nu_0}(\cdot)\label{tran2}\\
&&\ \ \ \ \ \ + \sum_{n=1}^{\infty}d^{\theta}_n(t)\int_{S^n}\nu^n(d\,x_1\times \cdots\times d\,x_n)
\Pi_{n+\theta, \frac{n}{\theta+n} \eta_n +\frac{\theta}{\theta+n}\nu_0}(\cdot),\nn
\eeqn
and $\Pi_{\theta,\nu_0}$ is the reversible measure of the process.
\end{theorem}

It is clear from this theorem that for any $t>0, $ $Y_t$ and $X_t$ are mixtures of gamma processes and the Dirichlet processes respectively.  A comparison between $Q_1(t,\mu,\cdot)$ and $Q_2(t,\nu,\cdot)$ reveals a termwise Gamma-Dirichlet structure between
\[
\Gamma^{C(-\la,t)}_{n+\theta, \frac{n}{\theta+n} \eta_n +\frac{\theta}{\theta+n}\nu_0}(\cdot)\ \mb{and}\ \Pi_{n+\theta, \frac{n}{\theta+n} \eta_n +\frac{\theta}{\theta+n}\nu_0}(\cdot).
\]

\vspace{0.3cm}
Let $C([0,\infty), M(S))$ denote the space of all $M(S)$-valued continuous functions on $[0,\infty)$ equipped with the topology of uniform convergence on compact sets. Define
\[
\varrho(\mu(\cdot))=\inf\{t>0: \mu_t(S)=0 \}, \ \mu_{\cdot} \in C([0,\infty), M(S))
\]
and
\[
C_0([0,\infty), M(S))=\{\mu_{\cdot}\in C([0,\infty), M(S)): \int_{0}^{\varrho(\mu(\cdot))}\frac{d\,u}{\mu_u}=\infty\}.
\]
\vspace{0.5cm}
Define
\[
\varsigma(t): C_0([0,\infty), M(S))\ra [0,\infty),\ t =\int_0^{\varsigma(t)}\frac{d\,u}{\mu_u(S)}.
\]
\vspace{0.3cm}
It is shown in \cite{Shiga90} that
the process $\frac{Y_{\varsigma(t)}}{Y_{\varsigma(t)}(S)}$ is the FVP starting at $\nu=\frac{\mu}{\mu(S)}$. A natural question raised in \cite{EtGri93b} is whether one can derive \rf{tran2} directly from \rf{tran1} using Shiga's normalization and random time change. Motivated by this problem, we obtain a direct derivation of $\{d^{\theta}_n(t): n =0,1,\ldots\}$ from $\{q_n^{\mu(S),\la}(t): n=0,1,\ldots\}$ through a random time change. This result is then used to link the distribution of $Y_t$ and $X_t$ at every fixed time $t >0$.

For notational simplicity, we write $N(t)$ for the time inhomogeneous Markov chain $N_{\mu(S)}(t)$. Define $\tau_t: [0,\infty)\rightarrow [0,\infty)$ by

\[
t=\int_{0}^{\tau_t}\frac{d\,u}{(N(u)\vee 1 +\theta-1)C(-\la,u)}.
\]
\vspace{0.3cm}

\begin{theorem}\label{ran-time-ch}
The process $N(\tau_t)$ is the embedded chain of Kingman's coalescent $D_t^{\theta}$.
\end{theorem}

\proof  Let $C(\hat{Z}_+)$ be the set of all continuous functions on $\hat{Z_+}$ and
\[
C_0 = \{f \in C(\hat{Z}_+): \lim_{n\ra \infty}n^2 (f(n-1)-f(n)) \ \mb{exists}\}.
\]

For any $t > 0$ and $f$ in $C_0$, define
\[
\Omega_t f(m)= \left\{\begin{array}{lr}0, & m=\infty \\ \frac{m}{2C(-\la,t)}(f(m-1)-f(m)),& m \in Z_+ \end{array}\right.
\]
\vspace{0.3cm}
Then $\Omega_t$ is the time dependent infinitesimal generator of $N(t)$. Hence for any $f$ in $C_0$ and $t,s >0$
\be\label{mar1}
f(N(t+s))-f(N(s))-\int_0^{t} \Omega_{s+u} f(N(s+u))d\,u
\ee
is a martingale with respect to the natural filtration of $\{N(t+s): t\geq 0\}$. Let $\tau_t(s)$ be given by
\[
t=\int_s^{s+\tau_t(s)}\frac{d\,u}{(N(u)\vee 1+\theta-1)C(-\la,u)}.
\]
Clearly $\tau_t(s)$ converges to $\tau_t$ as $s$ tends to zero. It follows from \rf{mar1} that
\beq
&& f(N(s+\tau_t(s)))-\int_0^{\tau_t(s)}\Omega_{s+u} f(N(s+u))d\,u\\
&&= f(N(s+\tau_t))\\
&&\ \ \ \ \ \ -\int_0^{t}(\Omega_{s+\tau_u(s)} f(N(s+\tau_u(s)))) (N(s+\tau_u(s))\vee 1 +\theta-1)C(-\la,s+\tau_u(s))d\,u\\
&& = f(N(s+\tau_t(s)))-\int_0^{t}\la_{N(s+\tau_u(s))}[f(N(s+\tau_u(s))-1)-f(N(s+\tau_u(s)))]d\,u
\eeq
is a martingale with respect to the filtration generated by $\{N(\tau_t(s)+s): t\geq 0\}$. Set
$$
\tilde{\Omega} f(m)=\left\{\begin{array}{cr}
                       \la_m (f(m-1)-f(m)) &  m \in Z_+,\\
                      \lim_{n\ra \infty}\la_n [f(n-1)-f(n)] & m=\infty
                     \end{array}
\right.
$$
It then follows by letting $s$ tend to zero that the process $N(\tau_t)$ is a solution to the martingale problem associated with the generator $\tilde{\Omega}$. On the other hand, based on the  proof of lemma 2.5 in \cite{EtGri93a}, the process $D_t^{\theta}$ is the unique solution of the martingale problem associated with $\tilde{\Omega}$. Therefore the theorem holds.

\hfill $\Box$

As an application of this result, we get the following derivation of the fixed time distribution of the FVP process from the MBI process.

\begin{theorem}\label{TFV}
For any $t>0$, let
\[
\nu(\theta,t)=\frac{N(t)\eta_{N(t)}+\theta\nu_0}{N(t)+\theta}.
\]
Then we have
\begin{eqnarray}
Y_t&\overset{d}{=}&{\cal Y}^{C(-\la,t)}_{N(t)+\theta,\nu(\theta,t)} \label{dy-GD1}\\
X_t &\overset{d}{=}&\frac{{\cal Y}^{C(-\la,t)}_{N(\tau_t)+\theta, \nu(\theta, \tau_t)}}{{\cal Y}^{C(-\la,t)}_{N(\tau_t)+\theta,\nu(\theta,\tau_t)}(S)}.\label{dy-GD2}
\end{eqnarray}
\end{theorem}

\noindent {\bf Remarks} In \cite{Shiga90}, the MBI process $Y_t$ can be represented as the sum of two independent Poisson clusters with one corresponding to the mass distribution of the descendants of the original population and the other the mass distribution of all immigrants.  The random time change in \cite{Shiga90} is for the whole process $Y_t$ while the random time change here is only for the Poisson number of descendants. The scaling parameter $C(-\la,t)$ plays no role in \rf{dy-GD2} and therefore the random time change involves $N(t)$ only.

\subsection{Reversibility}

 The Fleming-Viot process is a large class of probability-valued processes that describe the evolution of a population under the influence of mutation, selection, recombination, and random sampling. If there is no selection and recombination, then the FVP process $X_t$ is shown in \cite{lsy99} to be the only reversible Fleming-Viot process. The corresponding reversible measure is $\Pi_{\theta,\nu_0}$.
In \cite{Handa02}, the reversibility of $X_t$ is shown to be equivalent to
 the quasi-invariance of the Dirichlet process $\Pi_{\theta,\nu_0}$. In this subsection we study the reversibility of a class of branching diffusions with immigration and investigate the corresponding relation between reversibility and quasi-invariance.

For any $\mu_0$ in $M(S)$, $a(\cdot), b(\cdot)$ in $C(S)$ satisfying $a(\cdot)>0$ and $b(\cdot)>0$, consider the following generator

\be
\tilde{{\cal L}}F(\mu)=\int_S \mu(d\,s) a(s)\frac{\delta^2 F(\mu)}{\delta\mu(s)^2}
+\int_S (\mu_0(d\,s)-\mu(d\,s)b(s))\frac{\delta F(\mu)}{\delta \mu(s)}\label{Bran}
\ee
with domain
$$
\mathcal{D}(\tilde{{\cal L}})=\{\phi(\langle \mu,f_1\rangle,\ldots, \langle \mu,f_n\rangle);\mu \in M(S),\
n\geq1, f_i\in B(S), \phi\in C^2(\mathbb{R}^n)\}.
$$

It is known that the martingale problem associated with $\tilde{{\cal L}}$ is well-posed (cf. \cite{Shiga90}). The unique solution to the martingale problem is a diffusion process $Z_t$ with fixed time distribution characterized
by
\beqn
E_\mu[e^{-\langle Z_t,f\rangle}]&=&\exp\(-\left\langle \mu_0, a^{-1}\log\(1+\frac{a}{b}(1-e^{-bt})f\)\right\rangle \)\label{laplace}\\
&& \ \ \ \times \exp\( -\left\langle \mu, \frac{e^{-b t}}{1+a/b(1-e^{-bt})f}\right\rangle
\).\nn
\eeqn
\vspace{0.3cm}

The corresponding carr\'{e} du champ is
\beqn
 \Gamma(F, G)&=& \frac{1}{2}\{\tilde{{\cal L}}(FG)-G\tilde{{\cal L}}F-F \tilde{{\cal L}}G\}\label{cdc}\\
&=&\left\langle \mu, a(s)\frac{\delta F(\mu)}{\delta \mu(s)}\frac{\delta G(\mu)}{\delta \mu(s)}\right\rangle,\nn
\eeqn
 where $F,G\in\mathcal{D}(\tilde{{\cal L}})$. It follows from direct calculation that for any $F,G, H$ in $\mathcal{D}(\tilde{{\cal L}})$
\begin{equation}\label{gg}
\Gamma(FH, G)+\Gamma(HG, F)-\Gamma(H,FG)=2H\Gamma(F,G).
\end{equation}

\vspace{0.3cm}

For any $f$ in $B(S)$ and $\mu$ in $M(S)$, define
\[
S_f: M(S)\ra M(S),\  \mu(d\,s)\mapsto e^{f(s)}\mu(d\,s)
\]
and
\[
\Lambda(\mu,f)=\langle \mu_0,f/a\rangle-\langle \mu,(e^{f}-1)b/a\rangle.
\]

\vspace{0.3cm}
\begin{definition}
Let $\Xi$ be a probability on $M(S)$. The operator $\tilde{{\cal L}}$ is reversible with respect to $\Xi$ if
$$
\int F \tilde{{\cal L}}G \Xi(d\,\mu)=\int G \tilde{{\cal L}}F \Xi(d\,\mu),\qquad F,G\in \mathcal{D}(\tilde{{\cal L}}).
$$
The probability $\Xi$ is $\Lambda$-quasi-invariant with respect to the family of transformations $\{S_f: f \in B(S)\}$  if $\Xi$ and the image law $S_f(\Xi)$ of $\Xi$ under $S_f$ are mutually absolutely continuous with density given by
\[
\frac{d\,S_f(\Xi)}{d\,\Xi}(\mu)=e^{\Lambda(\mu,-f)}.
\]
\end{definition}

The next result shows that the reversibility and the $\Lambda$-quasi-invariance are equivalent.
\begin{theorem}\label{equivalent}
 A Borel probability measure  $\Xi$ on $\M(S)$ is reversible with respect to $\tilde{{\cal L}}$ if and only if $\Xi$ is $\Lambda$-quasi-invariant with respect to the family of transformations $\{S_f: f \in B(S)\}$.
\end{theorem}

\proof We first show that the reversibility implies the $\Lambda$-quasi-invariance. From the definition of the carr\'{e} du champ, it is clear that for any $F,G$ in ${\cal D}(\tilde{{\cal L}})$
\be\label{lap1}
\int_{M(S)}[\Gamma(F,G)+F\tilde{{\cal L}}G]\Xi(d\,\mu)=0.
\ee

For any $t \geq 0, g \in B(S)$, let $$\mu_t=S_{-atg}(\mu)$$ and
$$
F(\mu; f_1,\dots,f_n)=\phi(\langle \mu,f_1\rangle,\ldots,\langle \mu,f_n\rangle),$$
where $n \geq 1, \phi\in \mathcal{D}(\tilde{{\cal L}}), f_i\in B(S), i=1,2,\ldots$. It is clear that
\[
F(\mu_t; f_1,\dots,f_n)=F(\mu; e^{-a tg}f_1, \ldots, e^{-atg}f_n).
\]
\vspace{0.4cm}
By direct calculation
\be\label{qu-den5}
\frac{d\,F(\mu_t;f_1,\ldots,f_n)}{d\,t}= -\left\langle \mu, a g \frac{\delta F(\mu_t;f_1,\ldots,f_n)}{\delta \mu(s)}\right\rangle.
\ee
\vspace{0.3cm}
For $G(\mu)=\langle \mu,g\rangle$, we have
\be\label{qu-den6}
\tilde{{\cal L}}G(\mu)= \langle \mu_0, g \rangle -\langle \mu, b g \rangle
\ee
and
\beqn
\Lambda(\mu_t,a t g)&=& \langle \mu_0, g\rangle t- \langle \mu, \frac{b}{a}(1-e^{-atg})\rangle \label{qu-den7}\\
&=&\int_0^t \tilde{{\cal L}}G(\mu_u)d\,u\nn.
\eeqn

\vspace{0.3cm}

Set

\beq
F_t(\mu)&=&F(\mu_t; f_1,\ldots,f_n)e^{-\Lambda(\mu_t,atg)}\\
H(t)&=&\int_{M(S)}F_t(\mu)\Xi(d\,\mu).
\eeq
\vspace{0.3cm}
Since
\[
\frac{d\, \tilde{{\cal L}}G(\mu_t)}{d\,t}= -\left\langle \mu, a g\frac{\delta \tilde{{\cal L}}G(\mu_t)}{\delta \mu(s)}  \right\rangle,
\]
it follows from \rf{qu-den5} and \rf{qu-den7} that

\beq
\left\langle \mu, a g \frac{\delta F_t(\mu)}{\delta \mu(s)}\right\rangle &=& \left\langle \mu, \frac{\delta F(\mu_t; f_1,\ldots,f_n)}{\delta \mu(s)}\right\rangle e^{-\Lambda(\mu_t,atg)}- F_t(\mu)\int_0^t \left\langle\mu, a g\frac{\delta \tilde{{\cal L}}G(\mu_s)}{\delta \mu(s)}
\right\rangle d\,s\\
&=&-\frac{d\, F(\mu_t; f_1,\ldots,f_n)}{d\, t}e^{-\Lambda(\mu_t,atg)}+F_t(\mu)(\tilde{{\cal L}}G(\mu_t)-\tilde{{\cal L}}G(\mu)).
\eeq
\vspace{0.4cm}
Consequently,
\begin{eqnarray*}
H'(t)&=&\int_{M(S)} \left[\frac{d}{dt}F(\mu_t; f_1,\ldots,f_n)e^{-\Lambda(\mu_t,atg)}-F_t(\mu)\tilde{{\cal L}}G(\mu_t)
 \right]\Xi(d\mu)\\
&=&\int_{M(S)} \left[-\langle \mu,ag\frac{\delta F_t(\mu)}{\delta \mu(s)}\rangle-F_t(\mu)\tilde{{\cal L}}G(\mu)
\right]\Xi(d\mu)\\
&=&-\int_{M(S)} \left[\Gamma(F_t, G)+F_t(\mu)\tilde{{\cal L}}G(\mu)
\right]\Xi(d\mu)
\end{eqnarray*}
\vspace{0.3cm}

Since  $F_t, G\in\mathcal{D}(L)$, it follows from \rf{lap1} that
 $H'(t)=0$. In particular we have
$$
\int F(\mu_1;f_1,\ldots,f_n)e^{-\Lambda(\mu_1,a g)} \Xi(d\mu)=\int F(\mu_0;f_1,\ldots,f_n)\Xi(d\mu).
$$

\vspace{0.3cm}

Since $g$ is arbitrary and $a(x)>0$, we have for any $f\in B(S)$

$$
\int F(S_{-f}(\mu); f_1,\ldots,f_n)e^{-\Lambda(S_{-f}\mu, f)} \Xi(d\mu)=\int F(\mu;f_1,\ldots,f_n)\Xi(d\mu),
$$
\vspace{0.3cm}
 Replacing $F(\mu;f_1,\ldots,f_n)$ by $F(\mu; f_1,\ldots,f_n)e^{\Lambda(\mu,f)}$, we obtain
$$
\int_{M(S)} F(S_{-f}\mu;f_1,\ldots,f_n) \Xi(d\mu)=\int F(\mu; f_1,\ldots,f_n)e^{\Lambda(\mu,f)}\Xi(d\mu),
$$
which shows that $\Xi$ is $\Lambda$-quasi-invariant.

\smallskip

Next we show that the $\Lambda$-quasi-invariance implies reversibility.  Assume that $\Xi$ is $\Lambda$-quasi-invariant. Then we have $H'(t)=0$ for all $t\geq 0$. In particular, we have for any $g$ in $B(S)$ and $F$ in ${\cal D}(\tilde{{\cal L}})$
\beq
-H'(0)&=&\int_{M(S)}[\langle \mu, ag\frac{\delta F(\mu)}{\delta \mu(x)}\rangle+F(\mu)\tilde{{\cal L}}G(\mu)
]\Xi(d\mu)\\
&=&\int_{M(S)}[\Gamma(F,G)+F(\mu)\tilde{{\cal L}}G(\mu)
]\Xi(d\mu) =0,
\eeq
where as before $G(\mu)=\langle \mu, g\rangle$. For any $g_1,g_2$ in $B(S)$, let $G_i(\mu)=\langle\mu,g_i\rangle, i=1,2.$ Then it follows from \rf{gg} that
\beq
\int_{M(S)} F \tilde{{\cal L}}(G_1 G_2)\Xi(d\mu)&=&\int_{M(S)} F(2\Gamma(G_1,G_2)+G_2\tilde{{\cal L}} G_1+G_1 \tilde{{\cal L}} G_2)
\Xi(d\mu)\\
&=&\int_{M(S)} (\Gamma(F G_1, G_2)+\Gamma(F G_2, G_1)-\Gamma(F,G_1 G_2))\Xi(d\mu)\\
&&\ \ -\int_{M(S)} \Gamma(F G_1, G_2)+\Gamma(F G_2, G_1)\Xi(d\mu)\\
&=& -\int \Gamma(F,G_1 G_2)\Xi(d\mu).
\eeq
By induction,
\be\label{rev-final}
\int_{M(S)} F \tilde{{\cal L}}G \Xi(d\,\mu)=-\int_{M(S)}\Gamma(F,G)\Xi(d\,\mu)
\ee
for $G$ of the form $\prod_{i=1}^n \langle\mu, g_i\rangle, g_i \in B(S), n \geq 1$. Note that every function $\phi$ in $C^2(\mathbb{R}^n)$  can be approximated by polynomials under the Sobolev norm of order two, so the equality \rf{rev-final} holds for any $G$ in ${\cal D}(\tilde{{\cal L}})$. Consequently, $\Xi$ is reversible with respect to $\tilde{{\cal L}}$.

\hfill $\Box$

As an application of Theorem~\ref{equivalent}, we consider the reversibility of the process $Z_t$.
Letting $t$ tend to infinity in \rf{laplace}, it follows that the process $Z_t$ has a unique invariant distribution $\Gamma_{a^{-1}b,a^{-1}\mu_0}$ characterized by

\be\label{laplace1}
E^{\Gamma_{a^{-1}b,a^{-1}\mu_0}}[e^{-\langle \mu, h\rangle}]=\exp\(-\left\langle \mu_0, a^{-1}\log(1+\frac{a}{b}h)\right\rangle\),\ \ h \in B_+(S).
\ee
\vspace{0.3cm}

Next we establish the reversibility of
$Z_t$ by verifying the $\Lambda$-quasi-invariance of $\Gamma_{a^{-1}b,a^{-1}\mu_0}$.

\begin{theorem}\label{laplace2}
The law $\Gamma_{a^{-1}b,a^{-1}\mu_0}$ is $\Lambda$-quasi-invariant and therefore the process $Z_t$ is reversible with respect to $\Gamma_{a^{-1}b,a^{-1}\mu_0}$.
\end{theorem}

\proof For any $f$ in $B(S)$ and $h$ in $B_+(S)$, it follows from \rf{laplace1} that

\beq
&&E^{S_f(\Gamma_{a^{-1}b,a^{-1}\mu_0})}[e^{-\langle \mu, h\rangle}]\\
&&\ \ \ \ = E^{\Gamma_{a^{-1}b,a^{-1}\mu_0}}[e^{-\langle S_{f}(\mu), h\rangle}]\\
&&\ \ \ \ = E^{\Gamma_{a^{-1}b,a^{-1}\mu_0}}[e^{-\langle \mu, e^{f} h\rangle}]\\
&&\ \ \ \ = \exp\{-\left\langle \mu_0, a^{-1}\log(1+\frac{a}{b}h e^{f}) \right\rangle\}\\
&&\ \ \ \ =\exp\{-\langle \mu_0, a^{-1}f\rangle\}\exp\left(-\left\langle\mu_0, a^{-1}\log \(1+\frac{a}{b}(h +\frac{b}{a}(e^{-f}-1))\)\right\rangle\right)\\
&&\ \ \ \ =\exp\{-\langle \mu_0, a^{-1}f\rangle\}E^{\Gamma_{a^{-1}b,a^{-1}\mu_0}}[e^{-\langle \mu, (h +\frac{b}{a}(e^{-f}-1))\rangle}]\\
&&\ \ \ \ = E^{\Gamma_{a^{-1}b,a^{-1}\mu_0}}[e^{\Lambda(\mu,-f)}e^{-\langle \mu, h\rangle}],
\eeq
which yields the result.

\hfill $\Box$

\noindent {\bf Remarks} The $\Lambda$-quasi-invariance is established in \cite{TVY01} when $a$ and $b$ are constants. For the MBI process $Y_t$,  we have $a(s)=\frac{1}{2}$, $b(s)=\frac{\lambda}{2}$, and $\mu_0=\nu_0$.

\vspace{0.3cm}

We conclude this subsection with a derivation of the reversibility of the FVP process from the reversibility of the MBI process exploiting the Gamma-Dirichlet structure.

 \begin{theorem}
  The reversibility of the MBI process implies the reversibility of the FVP process.
 \end{theorem}
 \proof First recall that the domain of the generator ${\cal A}$ for the FVP process is given by
\[
{\cal D}({\cal A})=\{\phi(\langle \nu,f_1\rangle,\ldots, \langle \nu,f_n\rangle);\nu \in M_1(S),\
n\geq1, f_i\in B(S), \phi\in C^2(\mathbb{R}^n)\}.
\]
 For any $ \mu \in M(S)\setminus \{0\}$ and any $\Phi, \Psi$ in ${\cal D}({\cal A})$, define
 \beq
 r(\mu)&=&\langle\mu,1\rangle\\
 F(\mu)&=& r^3(\mu)\Phi(\frac{\mu}{\mu(S)})  \\
 G(\mu)&=& r^3(\mu) \Psi(\frac{\mu}{\mu(S)}),
 \eeq
 and $F(0)=G(0)=0$. Then it is clear that both $F$ and $G$ belong to ${\cal D}({\cal L})$. By the definition of ${\cal D}({\cal A})$, there are $n,m \geq 1, f_1,\ldots, f_n, g_1,\ldots,g_m \in B(S), \phi, \psi \in C^2(\mathbb{R}^n)$ such that
 \beq
 \Phi(\nu)&=&\phi(\langle\nu, f_1\rangle, \ldots, \langle\nu,f_n\rangle)\\
 \Psi(\nu)&=&\psi(\langle\nu, g_1\rangle, \ldots, \langle\nu,g_m\rangle).
 \eeq
\vspace{0.3cm}
 Fix a $\mu$ in $M(S)\setminus \{0\}$. As in section 2.2, we write $\hat{\mu} = \frac{\mu}{\mu(S)}$. Set
 \[
  \tilde{\Phi}(\mu)=\Phi(\hat{\mu}),\ \tilde{\Psi}(\mu)=\Psi(\hat{\mu}).
\]
\vspace{0.3cm}
Then by direct calculation
 \beq
 \frac{\delta F(\mu)}{\delta \mu(s)} &=& 3r^2(\mu)\tilde{\Phi}(\mu) +r^3(\mu)\frac{\delta \tilde{\Phi}(\mu)}{\delta \mu(x)} \\
 &=& 3r^2(\mu)\tilde{\Phi}(\mu) +r^2(\mu)\sum_{i=1}^n\partial_i \phi (f_i-\langle \hat{\mu}, f_i\rangle)
 \eeq
and
\beq
\frac{\delta^2 F(\mu)}{\delta \mu(s)^2} &=& 6r(\mu)\tilde{\Phi}(\mu) +4r(\mu)\sum_{i=1}^n\partial_i \phi (f_i-\langle \hat{\mu}, f_i\rangle)\\
&&+r(\mu)\sum_{i,j=1}^n \partial_{ij}\phi (f_i-\langle \hat{\mu}, f_i\rangle)(f_j-\langle \hat{\mu}, f_j\rangle).
 \eeq
\vspace{0.3cm}
Substituting this into \rf{gen-mbi}, yields

 \be\label{delta1}
 {\cal L}F(\mu) = r^2(\mu) {\cal A}\Phi(\hat{\mu})+ 3r^2(\mu)(\frac{\theta-\lambda r(\mu)}{2}+2)\Phi(\hat{\mu}).
  \ee
Similarly
 \be\label{delta2}
 {\cal L}G(\mu) = r^2(\mu) {\cal A}\Psi(\hat{\mu})+ 3r^2(\mu)(\frac{\theta-\lambda r(\mu)}{2}+2)\Psi(\hat{\mu}).
  \ee

\vspace{0.5cm}
For $\mu=0$, we have ${\cal L}F(\mu)={\cal L}G(\mu)=0$. By the reversibility of the MBI process

\beq
0&=&\int_{M (S)} [G(\mu){\cal L}F(\mu)- F(\mu){\cal L}G(\mu)]\Gamma^{\beta}_{\theta,\nu_0}(d\,\mu)\\
&=& \int_{M (S)} r^5(\mu) [\Psi(\hat{\mu}){\cal A}\Phi(\hat{\mu})-\Phi(\hat{\mu}){\cal A}\Psi(\hat{\mu})]\Gamma^{\beta}_{\theta,\nu_0}(d\,\mu),
\eeq
which, combined with Gamma-Dirichlet structure, implies
\beq
&&\int_{M (S)} r^5(\mu) [\Psi(\hat{\mu}){\cal A}\Phi(\hat{\mu})-\Phi(\hat{\mu}){\cal A}\Psi(\hat{\mu})]\Gamma^{\beta}_{\theta,\nu_0}(d\,\mu)\\
&&\ \ \ \ =\int_0^{\infty} \frac{1}{\Gamma(\theta)\beta^{\theta}}x^{\theta+4}e^{-x/\beta}d\,x \int_{M_1(S)}[\Psi(\hat{\mu}){\cal A}\Phi(\hat{\mu})-\Phi(\hat{\mu}){\cal A}\Psi(\hat{\mu})]\Pi_{\theta,\nu_0}(d\,\hat{\mu})\\
&&\ \ \ \ =\frac{\beta^5\Gamma(5+\theta)}{\Gamma(\theta)}\int_{M_1(S)}[\Psi(\hat{\mu}){\cal A}\Phi(\hat{\mu})-\Phi(\hat{\mu}){\cal A}\Psi(\hat{\mu})]\Pi_{\theta,\nu_0}(d\,\hat{\mu})=0.
\eeq

\vspace{0.5cm}
Therefore the FVP process is reversible with respect to $\Pi_{\theta,\nu_0}$.

\hfill $\Box$

\noindent {\bf Remark} It is not clear whether the reversibility of the MBI process follows from the reversibility of the FVP process.

      \end{document}